# The Minimum $L_2$-Distance Projection onto the Canonical Simplex: A Simple Algorithm

Hans J. H. Tuenter

We consider the problem of the minimum distance projection of a given point $a$ in n-dimensional space onto the canonical simplex.

This problem can be formulated as follows:

$$\min \sum_{i=1}^{n} (x_i - a_i)^2$$

$$\text{s.t.} \sum_{i=1}^{n} x_i = 1,$$

$$x_i \geq 0, \ \forall i.$$

This is a problem that occurs in the setting of credit risk, where one deals with stochastic matrices that describe transition probabilities between different credit ratings, and where one wants to determine the roots of these matrices, or close approximations to them. A recent *ARQ* article by Kreinin and Sidelnikova (2001) gives the genesis and a more complete description of this problem.

## Solution

Without loss of generality we may assume that $a_1 \geq a_2 \geq \ldots \geq a_n$, since this is just a matter of reordering the data. If $x^*$ denotes the optimum solution, then the following lemmata allow us to reduce the problem to a univariate one.

**Lemma 1**. There is an index $m$ such that $x_i^* > 0$, for $1 \leq i \leq m$, and $x_i^* = 0$, for $i > m$.

**Proof**. Suppose that we have a feasible solution $x$ with $x_i = 0$ and $x_{i+1} > 0$, for some index $i$.

Now, consider the solution $\tilde{x}$, constructed from $x$ by switching the $i$th and $(i+1)$th components. This solution is feasible and has distance to the canonical simplex of

$$d(\tilde{x}) = 2(a_{i+1} - a_i)x_{i+1} + d(x).$$

But, as $a_i \geq a_{i+1}$, we have constructed a solution with distance less then or equal to that of the solution we started out with. This shows that we can limit ourselves to feasible solutions where a positive entry is never preceded by a zero entry, and thus proves the lemma.

**Lemma 2**. All positive entries of $x^*$ are of the form $x_i^* = a_i + \lambda$.

**Proof**. If there is only one positive element in $x^*$, the lemma is trivial, so assume that there are two elements $x_i^*$ and $x_j^*$ that are positive. Now consider the solution $\tilde{x}$, constructed from $x$ by perturbing the $i$th and $j$th components:

$\tilde{x}_i = x_i^* - \varepsilon$ and $\tilde{x}_j = x_j^* + \varepsilon$. For $\varepsilon$ sufficiently close to zero, this gives a feasible solution with distance to the canonical simplex of

$$d(\tilde{x}) = d(x) + 2\varepsilon^2 + 2\varepsilon(x_j^* - x_i^* - a_j + a_i).$$

If $x_j^* - x_i^* - a_j + a_i \neq 0$ holds, then one can always choose $\varepsilon$ such that $\tilde{x}$ has a lesser distance to the canonical simplex than $x^*$. But, this contradicts the optimality of $x^*$, and, thus, one must have $x_j^* - a_j = x_i^* - a_i$, which proves the lemma.

Combining the lemmata, one can reformulate the problem as an optimization problem over $\lambda$ and $m$.

$$\min \ m\lambda^2 + \sum_{i=m+1}^{n} a_i^2$$

$$\text{s.t.} \ m\lambda + \sum_{i=1}^{m} a_i = 1$$

$$\lambda \geq -a_m,$$
$$m \in \{1, 2, \ldots, n\}.$$

Substituting for $\lambda$ reduces it to a problem over the index $m$ only:

**Technical notes**

$$\min \quad \frac{1}{m}\left(1 - \sum_{i=1}^{m} a_i\right)^2 + \sum_{i=m+1}^{n} a_i^2$$

$$\text{s.t.} \quad \sum_{i=1}^{m} (a_i - a_m) \leq 1,$$

$$m \in \{1, 2, \ldots, n\}.$$

Now, consider the sequence

$$S_m = \sum_{i=1}^{m} (a_i - a_m).$$

This sequence is non-decreasing, as can be seen by formulating it in terms of the recursion

$$S_m = S_{m-1} + (m-1)(a_{m-1} - a_m),$$

with initial condition $S_1 = 0$. Denoting the objective function by $f(m)$, one can verify that it satisfies the recursion

$$f(m) = f(m-1) - \frac{1}{m(m-1)}[1 - S_m]^2,$$

with initial condition $f(1) = 1 - 2a_1 + \Sigma_{i=1}^{n} a_i^2$, and one immediately sees that the sequence $f(m)$ is non-increasing. These observations allow us to solve the problem by determining $m^*$ as the largest index $m$ that satisfies $S_m \leq 1$. Such an $m^*$ exists, since $S_1 = 0$, and is easily determined using the recursion for $S_m$. This then gives

$$\lambda^* = \frac{1}{m^*}\left(1 - \sum_{i=1}^{m^*} a_i\right) = \frac{1 - S_{m^*}}{m^*} - a_{m^*},$$

and the optimum solution as

$$x_i^* = \begin{cases} a_i + \frac{1}{m^*}\left(1 - \sum_{i=1}^{m^*} a_i\right), & \text{for } 1 \leq i \leq m^*, \\ 0, & \text{for } i > m^*. \end{cases}$$

Note that the optimum value for $m$ need not be unique, as $S_{m^*} = 1$ may hold.

## Illustrative example

To illustrate the algorithm, we use data from the matrix given in Table 3 of the article by Kreinin and Sidelnikova (2001). This matrix represents the square root of an annual credit-rating transition matrix, and contains negative elements that one would like to eliminate in order to obtain a six-month transition matrix. We take its first row and determine the permutation $\pi$ that orders the elements in descending order to obtain the vector:

$$\boldsymbol{a} = [0.947127, 0.051650, 0.001145,$$
$$0.000140, 0.000000, -0.000005,$$
$$-0.000006, -0.000050].$$

Applying the recursion to calculate $S_m$ we obtain the sequence:

$$0, 0.895477, 0.996487, 0.999502,$$
$$1.00062, 1.000087, 1.000093, 1.000401,$$

so that $m^* = 4$. We can now calculate $\lambda^* = -1.55 \times 10^{-5}$, and determine the optimum projection as

$$\boldsymbol{x}^* = [0.9471115, 0.0516345, 0.0011295,$$
$$0.0001245, 0, 0, 0, 0],$$

with a distance of $3.522 \times 10^{-9}$. Applying the inverse permutation $\pi^{-1}$ to $\boldsymbol{x}^*$ renders the first row of the regularized transition matrix. Applying the algorithm to the other rows of the matrix, one obtains the six-month transition matrix, as given in Table 4 of Kreinin and Sidelnikova (2001).

## Conclusions

In this Technical Note, we have considered and solved the problem of determining the minimum distance projection in the $L_2$-norm from an arbitrary point in an $n$–dimensional, Euclidian space onto the canonical simplex. It is shown that this problem reduces to a univariate problem that can be solved by a simple algorithm.

Finally, it should be remarked that the lemmata can also be derived through an application of the celebrated Kuhn-Tucker conditions (see Kuhn and Tucker 1951), but that we prefer to present a derivation from first principles.